\documentclass[11pt,english]{amsart}
\usepackage[a4paper]{geometry}
\usepackage{amssymb}
\usepackage{amsmath}
\usepackage{fancyhdr}
\usepackage{color}
\usepackage{enumerate}
\newtheorem{theorem}{Theorem}[section] 
\newtheorem{lemma}[theorem]{Lemma}     
\newtheorem{corollary}[theorem]{Corollary}
\newtheorem{proposition}[theorem]{Proposition}

\renewcommand{\Re}{\operatorname{Re}}
\renewcommand{\Im}{\operatorname{Im}}

\newtheorem{remark}[theorem]{Remark}
\newtheorem{remarks}[theorem]{Remarks}

\newtheorem*{theorem*}{Theorem}

\usepackage[bookmarksopen, naturalnames]{hyperref}

\makeatother

\begin{document}

\title[Universality of general Dirichlet series]{Universality of general Dirichlet series with respect to partial sums}
\author{A. Mouze}
\address{Augustin Mouze, \'Ecole Centrale de Lille, Univ. Lille, CNRS, UMR 8524 - Laboratoire Paul Painlev\'e  F-59000 Lille, France}
\email{augustin.mouze@univ-lille.fr}


\begin{abstract} 
	This paper investigates the universal approximation properties of general Dirichlet series $\sum_{n\geq 1} a_n e^{-\lambda_n s}$ absolutely convergent in the half-plane $\{\Re(s)>0\}$ under the condition $\lambda_{n+1}-\lambda_n \to 0$. We establish the existence of such series whose partial sums exhibit universal approximation on the closed left half-plane $\{\Re(s)\leq 0\}$. Furthermore, we prove that the Riesz typical means of these universal series preserve the universal approximation property, which gives, in particular, the automatic universality of logarithmic means of universal ordinary Dirichlet series. We also construct general Dirichlet series absolutely convergent in $\{\Re(s)>0\}$ that possess the universal approximation property on a restricted strip $\{\Re(s)\leq 0\}$ while their partial sums diverge to infinity for almost every $s$ in $\{\Re(s)<\sigma\}$.
\end{abstract}

\maketitle

\section{Introduction} A general Dirichlet series is a function of the form
\[
f(s)=\sum_{n\geq 1}a_ne^{-\lambda_n s},
\] 
where $(a_n)\subset\mathbb{C}^{\mathbb{N}}$ and $\lambda:=(\lambda_n)$ is an increasing sequence of nonnegative real numbers tending to $+\infty$ called a \textit{frequency}. The case $(\lambda_n)=(\log(n))$ gives the natural example of ordinary Dirichlet series. To any Dirichlet series 
$f=\sum_{n\geq 1}a_n e^{-\lambda_n s}$, one can associate its abscissa of absolute convergence 
\[
\sigma_{a}(f)=\inf\{\sigma\in\mathbb{R}:\sum_{n\geq 1}\vert a_n\vert e^{-\lambda_n \sigma}\hbox{ converges}\},
\]
and the sequence of its partial sums $(D_n(f))$, i.e. $D_n(f)(s)=\sum_{k=1}^na_k e^{-\lambda_k s}.$ Let $\mathbb{C}_+:=\{s\in\mathbb{C}:\Re(s)>0\}$ (resp. $\mathbb{C}_-:=\{s\in\mathbb{C}:\Re(s)<0\}$) be the right (resp. left) half-plane. Let us consider the Fr\'echet space $\mathcal{D}_a^{\lambda}(\mathbb{C}_+)$ of 
absolutely convergent Dirichlet series $\sum_{n\geq 1}a_ne^{-\lambda_n s}$ in $\mathbb{C}_+$ endowed with its natural topology given 
by the semi-norms $\Vert \sum_{n\geq 1}a_n e^{-\lambda_n s}\Vert_{\sigma}=\sum_{n\geq 1}\vert a_n\vert e^{-\lambda_n \sigma},$ 
$\sigma>0.$ When $\lambda=(\log(n))$, we will denote  $\mathcal{D}_a^{\lambda}(\mathbb{C}_+)$ by  $\mathcal{D}_a(\mathbb{C}_+)$. In this paper, we are interested in universal approximation properties of general Dirichlet series with respect to partial sums. 

For a compact set $K\subset\mathbb{C},$ $A(K)$ will stand for the space of all continuous 
functions on $K$ which are holomorphic in its interior. The notion of universal Taylor series was independently introduced by Luh 
\cite{Luh} and Chui and Parnes \cite{ChuiPar}, and was strengthened by Nestoridis \cite{Nes}, who showed the existence of 
holomorphic function $f$ on the unit disc $\{z\in\mathbb{C}:\vert z\vert <1\}$ whose Taylor series $\sum_{n\geq 0}a_nz^n$ at $0$ 
satisfies the following universal approximation property: for every compact set $K\subset\{z\in\mathbb{C};\vert z\vert\geq 1\}$ 
with connected complement and for every function $h\in A(K),$ there exists 
a subsequence $(\mu_n)\subset\mathbb{N}$ such that 
$\sum_{k=0}^{\mu_n}a_kz^k$ converges to $h$ uniformly on $K,$ as $n$ tends to infinity. Thus the sequence of partial sums of the Taylor series at $0$ of a holomorphic function on the unit disc can have the worst possible behavior, even on the unit circle.
Since Nestoridis' theorem, many results on universal Taylor series in the complex plane 
appeared. We refer the reader to \cite{bgnp} and the references therein. In 2005, Bayart established the existence of \textit{universal ordinary Dirichlet series} in $\mathcal{D}_a(\mathbb{C}_+)$ in the sense of Nestoridis: for any ``admissible'' compact set $K\subset \overline{\mathbb{C}_{-}}$ with connected complement and for any function $h\in A(K),$ there exists an increasing sequence $(\mu_n)\subset\mathbb{N}$ such that $\sup_{s\in K}\vert D_{\mu_n}(f)(s)-h(s)\vert \rightarrow 0$ as $n$ tends to infinity \cite{B1}. In this context, the condition for a compact set to be admissible is a geometric restriction on the width of its connected components. In \cite{ABGM} the authors relaxed this assumption. Further results on universal Dirichlet series can be found in \cite{bgnp}, \cite{Mou1},\cite{Mou2} \cite{Mou3} or \cite {GardiMano}. Other types of universality for Dirichlet series had previously been established, particularly with respect to vertical translations for the Riemann zeta function \cite{Voro}. We refer the reader to the survey \cite{Matsu}. Most of these results concern ordinary Dirichlet series. Recently several new universality results with respect to vertical translations or rearrangements of general Dirichlet series have been established \cite{Bay3,BayK}.\\
 
From now on, we will say that a sequence of holomorphic function $(f_n)$ satisfies \textit{the universal approximation property on a subset $E\subset\mathbb{C}$} if: for every compact set $K\subset E$ 
with connected complement and for every function $h\in A(K),$ there exists a subsequence $(\mu_n)\subset\mathbb{N}$ such that the sequence $(f_{\mu_n}(z))$ converges to $h$ uniformly on $K,$ as $n$ tends to infinity.\\

In this paper, we introduce the notion of universal general Dirichlet series in the sense of Nestoridis provided that the frequency $(\lambda_n)$ satisfies the assumption $\lambda_{n+1}-\lambda_n\rightarrow 0$, as n goes to infinity. To do this, we extend the approximation property given by Theorem 3.1 of \cite{ABGM} to general Dirichlet series and we obtain the existence of general Dirichlet series in $\mathcal{D}_a^{\lambda}(\mathbb{C}_+)$ whose sequence of partial sums satisfies the universal approximation property on $\overline{\mathbb{C}_-}$. We refer the reader to Theorem \ref{density}.\\

Then we show that for a universal general Dirichlet series with frequency $(\lambda_n)$ its Riesz typical means $(R,\lambda_n,k)$ ($k\in\mathbb{N}$) possess the same universal approximation property. The definitions of the Riesz typical means of a general Dirichlet series are given in Section \ref{RTM}. Notice that in particular in the case of ordinary Dirichlet series this result gives an answer to an open problem on the automatic universality of logarithmic means of a universal ordinary Dirichlet series. Moreover it extends Theorem 5.2 of \cite{Mou2} where the author showed that a universal ordinary Dirichlet series cannot be $(R,\log(n),1)$ summable at any point of its line of convergence. Here the automatic universality of logarithmic means of universal Dirichlet series shows that this one cannot be $(R,\log(n),k)$ summable at any point of its line of convergence for all $k\in\mathbb{N}$.\\ 

Finally we are interested in the notion of restricted universality. We know that there exists a ordinary Dirichlet series whose sequence of partial sums satisfies the universal approximation property on a strip $\{z\in\mathbb{C}:\sigma < \Re(z)\leq 0\}$ ($\sigma<0$) and tends to infinity for almost every (in the sense of Lebesgue measure) $z$ in the left half-plane $\{z\in\mathbb{C}:\Re(z)<\sigma-\frac{1}{2}\}$. It can be viewed as a weak Dirichlet version of a theorem of Kahane and Melas for universal Taylor series, which asserts the existence of holomorphic function $\sum_{k\geq 0}a_kz^k$ on the unit disc such that the sequence of its partial sums satisfies the universal approximation property on the annulus $\{z\in\mathbb{C}:1\leq \vert z\vert \leq a\}$ ($a>1$), and tends to infinity almost everywhere on $\{z\in\mathbb{C}:\vert z\vert >a\}$ \cite{MelKah}. Here we establish there exist general Dirichlet series in $\mathcal{D}_a^{\lambda}(\mathbb{C}_+)$ which possess the universal approximation property on a strip $\{s\in\mathbb{C}:\sigma\leq \Re(s)\leq 0\}$ and whose partial sums tend almost everywhere to infinity on $\{s\in\mathbb{C}:\Re(s)<\sigma\}$.  This improves the existing results for ordinary Dirichlet series and answers a question from \cite{Mou3}. We refer the reader to Theorem \ref{thmrestricted}. 
\vskip4mm

Let us note that all the results of this paper hold for general Dirichlet series whose frequency $(\lambda_n)$ satisfies $\lambda_{n+1} - \lambda_n\rightarrow 0$ as $n$ tends to infinity. For example, the statements apply to the following cases:
\begin{enumerate}[(a)]
	\item $\lambda_n= \log(n)$, i.e. the case of ordinary Dirichlet series;
	\item $\lambda_n=\log(p_n)$, where $(p_n)$ is the sequence of prime numbers, i.e. the case of ordinary Dirichlet series whose support is the set of prime numbers;
	\item $\lambda_n=n^a$, with $0<a<1$.
	\end{enumerate}
Let us also observe that there cannot exist a general Dirichlet series in $\mathcal{D}_a^{\lambda}(\mathbb{C}_+)$ with $\lambda=(n)$ satisfying the universal approximation property on $\overline{\mathbb{C}_-}$. Indeed, since any Dirichlet polynomial $P=\sum_{k=1}^n a_k e^{-ks}$ satisfies $P(i\pi)=P(-i\pi)$, the partial sums of such a general Dirichlet series cannot approximate, on a compact set containing both $i\pi$ and $-i\pi$, a function $g$ taking two different values at these points. \\

The paper is organized as follows. In Section \ref{gdsuniv}, we establish the existence of universal general Dirichlet series in $\mathcal{D}_a^{\lambda}(\mathbb{C}_+)$ with respect to partial sums, under the assumption $\lambda_{n+1}-\lambda_n\rightarrow 0$. Then Section \ref{properties} is divided into two parts: Subsection \ref{RTM} focuses on the automatic universality of Riesz typical means for such series, while Subsection \ref{restruniv} addresses restricted universality. 

\section{Universal general Dirichlet series}\label{gdsuniv} In this section, we extend Theorem 3.6 of \cite{ABGM}, which shows the existence of ordinary Dirichlet series in $\mathcal{D}_a(\mathbb{C}_+)$ whose sequence of partial sums satisfies the universal approximation property on $\overline{\mathbb{C}_-}$, to general Dirichlet series. The key ingredient for the construction of such universal ordinary Dirichlet series is a simultaneous polynomial approximation theorem \cite[Theorem 3.1]{ABGM}. Moreover, the existence of such a simultaneous polynomial approximation result is a necessary and sufficient condition for the existence of universal series in general \cite{bgnp}. Here we show that Theorem 3.1 of \cite{ABGM} remains true in the case of general Dirichlet series provided that the frequency $(\lambda_n)$ satisfies the condition $\lambda_{n+1}-\lambda_n\rightarrow 0$ as $n\to \infty$. The proof of this result, inspired by that of \cite[Theorem 3.1]{ABGM}, needs some preparation. We first give two lemmas on entire functions of exponential type.

\begin{lemma}\label{lemmaLexpo0} Let $h$ be an entire function of exponential type. Let $N\in \mathbb{N}$, $x>0$ and $(t_i(x))_{0\leq i\leq N}$ be an increasing sequence of $N+1$ positive real numbers with $x\in [t_0(x),t_N(x)]$ and $t_N(x)-t_0(x)\leq x/2$. Then, there exist $\alpha,A>0$ such that
	\[
	\vert h(x)\vert\leq \left(\sum_{i=0}^{N}\vert L_i(x)\vert\right)\left(\sup_{0\leq i\leq N}\vert h(t_i(x))\vert + 2Ae^{2\alpha x}\left(\frac{t_N(x)-t_0(x)}{x}\right)^{N+1}\right),
	\]
where $(L_i)_{0\leq i\leq N}$	is the family of Lagrange interpolation polynomials associated to the sequence $(t_i(x))_{0\leq i\leq N}$.
\end{lemma}
\begin{proof} By hypothesis, there exist $A,\alpha>0$ such that 
	\begin{equation}\label{expotype0}
		\forall z\in\mathbb{C},\quad\vert h(z)\vert\leq Ae^{\alpha\vert z\vert}.
	\end{equation}
Since $h$ is an entire function, we can write on $\vert z-x\vert\leq x$
	\[
	\begin{array}{rcl}\displaystyle h(z)=\sum_{k\geq 0}a_k(z-x)^k&=&\displaystyle\sum_{k=0}^{N}a_k(z-x)^k+\sum_{k\geq N+1}a_k(z-x)^k\\&:=&\displaystyle P_N(z)+\sum_{k\geq N+1}a_k(z-x)^k,\end{array}
	\]	
	with, by Cauchy's formula, taking into account (\ref{expotype0}) and the fact that the inequality $\vert z-x\vert\leq x$ implies $\vert z\vert\leq 2x$, 
	\[
	\forall k\geq 0,\quad \vert a_k\vert \leq\frac{Ae^{2\alpha x}}{x^k}
	\]
Now for all complex number $z$ such that $\vert z-x\vert\leq (t_N(x)-t_0(x))$ the following estimate holds
\begin{equation}\label{estimsg}
	\begin{array}{rcl}\displaystyle\vert P_N(z)-h(z)\vert&\leq&\displaystyle \sum_{k\geq N+1}Ae^{2\alpha x}\left(\frac{t_N(x)-t_0(x)}{x}\right)^k\\&\leq&\displaystyle Ae^{2\alpha x}\left(\frac{t_N(x)-t_0(x)}{x}\right)^{N+1}\frac{1}{1-\frac{t_N(x)-t_0(x)}{x}}\\
		&\leq& \displaystyle 2Ae^{2\alpha x}\left(\frac{t_N(x)-t_0(x)}{x}\right)^{N+1},\end{array}
\end{equation}
where the last inequality comes from the fact that $t_N(x)-t_0(x)\leq x/2$. Moreover observe that $h(x)=P_{N}(x)$. Given the set of $N+1$ reals $(t_{i}(x))_{0\leq i\leq N}$, let us consider the associated family of Lagrange polynomials $(L_i)_{0\leq i\leq N}$ of degree $N$ defined by 
	\[
	\forall i=0,\dots, N,\quad L_i(z)=\frac{\prod\limits_{\substack{j=0 \\ j\ne i}}^{N}(z-t_j(x))}{\prod\limits_{\substack{j=0 \\ j\ne i}}^{N}(t_i(x)-t_j(x))}.
	\]
	Therefore the following decomposition holds
	\[
	P_{N}(z)=\sum_{i=0}^{N}P_{N}(t_i(x))L_i(z).
	\]
	Note that for all $i=0,\dots, N$, $\vert t_i(x)-x\vert\leq t_N(x)-t_0(x)$. Combining this observation with (\ref{estimsg}), we get
	\begin{equation}\label{eqlag01}
		\begin{array}{rcl}\displaystyle\vert h(x)\vert=\vert P_{N}(x)\vert&\leq& \displaystyle\sum_{i=0}^{N}\vert P_{N}(t_i(x))\vert \vert L_i(x)\vert\\ &\leq& \displaystyle\sum_{i=0}^{N}\left(\vert P_{N}(t_i(x))-h(t_i(x))\vert +\vert h(t_i(x))\vert\right)\vert L_i(x)\vert\\&\leq&\displaystyle\left(\sum_{i=0}^{N}\vert L_i(x)\vert\right)\left(\sup_{0\leq i\leq N}\vert h(t_i(x))\vert+ 2Ae^{2\alpha x}\left(\frac{t_N(x)-t_0(x)}{x}\right)^{N+1} \right)
		\end{array}
	\end{equation}
Th proof is complete. 
\end{proof}

From the previous lemma, we derive the following result on entire functions of exponential type whose asymptotic behavior is known along an increasing unbounded sequence $(\lambda_n)$ satisfying the condition $\lambda_{n+1}-\lambda_n\rightarrow 0$ as $n\to\infty$. 

\begin{lemma}\label{lemmaLexpo} Let $(\lambda_n)$ be a frequency such that $(\lambda_{n+1}-\lambda_n)$ tends to $0$. Let $\mathcal{L}$ be an entire function of exponential type such that $\limsup\limits_{n\rightarrow\infty}\vert \mathcal{L}(\lambda_n)\vert e^{\sigma\lambda_n}<+\infty$. Then
	\[
	\limsup_{x\rightarrow +\infty}\frac{\log\vert\mathcal{L}(x)\vert}{x}\leq -\sigma.
	\] 	
\end{lemma}

\begin{proof} By hypothesis, there exist $A,\alpha>0$ such that 
	\begin{equation}\label{expotype}
		\forall z\in\mathbb{C},\quad\vert \mathcal{L}(z)\vert\leq Ae^{\alpha\vert z\vert}.
	\end{equation}
	Let $x>0$ sufficiently large. We need to estimate the term $\vert\mathcal{L}(x)\vert$ from the values of $\mathcal{L}$ along the sequence $(\lambda_{n})$. We will exploit the fact that $\lambda_{n+1}-\lambda_n\rightarrow 0$. We proceed by interpolation over a well-chosen small interval of length less than $x/2$ containing $x$ and $N_x$ terms of the sequence, where $\lim_{x\to\infty}N_x/x>0$. To this end, we will apply Lemma \ref{lemmaLexpo0}. First, we need some preparation. Let $\eta$ be in $(0,\sigma)$ and $c\in (0,1/4)$ such that 
	\begin{equation}\label{prepa1}
		c\log (4)\leq \eta.
	\end{equation} 
	By hypothesis, there exists a positive integer $n_0$ such that 
	\begin{equation}\label{prepa0}
	\hbox{for all }n\geq n_0,\quad \lambda_{n}> 1 \hbox{ and } \lambda_{n+1}-\lambda_n< 1.
	\end{equation} 
For $x >0 $, let us define
	\[
g(x)=\sup(\lambda_{k+1}-\lambda_k: \lambda_k\geq x/2).
\]
Since $(\lambda_n)$ is a frequency such that $(\lambda_{n+1}-\lambda_n)$ tends to $0$, we have $g(x)\rightarrow 0$ as $x\to\infty$.	Fo all $x>0$ sufficiently large, we can assume that 
\begin{equation}\label{equprelim}
x>\lambda_{n_0},\quad  g(x)\leq 1,\quad cx>1,\quad -c\log(cg(x))\geq 2\alpha+\sigma.
\end{equation}	
We set 
	\begin{equation}\label{prepa1bis}
	N_x:=\lfloor cx\rfloor\quad \hbox{ and }\quad \lambda_{m_0(x)}=\max(\lambda_k: \lambda_k\leq x).
	\end{equation}
Observe that $m_0(x)\geq n_0$. Therefore the following inequalities hold	
	\begin{equation}\label{prepa2}
		\lambda_{m_0(x)}\leq x<\lambda_{m_0(x)+1}\leq\lambda_{m_0(x)}+1<2\lambda_{m_0(x)},
	\end{equation}
which implies
\begin{equation}\label{prepa22}
	\lambda_{m_0(x)}\geq x/2.
\end{equation}
Let us define
\begin{equation}\label{prepa2300}
	\lambda_{m_1(x)}=\min(\lambda_k: \lambda_k\geq \lambda_{m_{0}(x)}+g(x)).
\end{equation}
We get 
\[
m_0(x)+1\leq m_1(x),\quad \lambda_{m_{0}(x)}+g(x)\leq \lambda_{m_1(x)}\quad\hbox{ and }\quad \frac{x}{2}\leq\lambda_{m_0(x)}\leq \lambda_{m_1(x)-1}<\lambda_{m_0(x)}+g(x).
\]
Therefore we have 
\[
\lambda_{m_{0}(x)}+g(x)\leq\lambda_{m_1(x)}=\lambda_{m_1(x)-1}+(\lambda_{m_1(x)}-\lambda_{m_1(x)-1})< \lambda_{m_0(x)}+2g(x)
\]
Now we proceed by induction to obtain an increasing sequence $(m_j(x))_{j=0,\dots,N_x}$ of positive integers such that for $j=1,\dots,N_x$, 
\begin{equation}\label{prepa231}
	\lambda_{m_{j-1}(x)}+g(x)\leq \lambda_{m_j(x)}<\lambda_{m_{j-1}(x)}+2g(x). 
\end{equation}
Indeed suppose that $m_{j-1}(x)$ has been constructed. It suffices to define 
\begin{equation}\label{prepa2301}
	\lambda_{m_j(x)}=\min(\lambda_k:\lambda_k\geq \lambda_{m_{j-1}(x)}+g(x)). 
\end{equation}
By definition, we have 
\[
m_{j-1}(x)+1\leq m_j(x),\quad \lambda_{m_{j-1}(x)}+g(x)\leq \lambda_{m_j(x)},\quad \lambda_{m_j(x)-1}<\lambda_{m_{j-1}(x)}+g(x),
\]
which implies $x/2\leq\lambda_{m_0(x)}<\dots <\lambda_{m_{j-1}(x)}<\lambda_{m_{j}(x)}$ and 
\[
\lambda_{m_{j-1}(x)}+g(x)\leq\lambda_{m_j(x)}=\lambda_{m_j(x)-1}+(\lambda_{m_j(x)}-\lambda_{m_j(x)-1})< \lambda_{m_{j-1}(x)}+2g(x).
\]
Thus (\ref{prepa231}) is proved. By construction we have 
	\[
	x\in [\lambda_{m_0(x)},\lambda_{m_1(x)})
	\]
and, by the property (\ref{prepa231}), for all $i,j=0,\dots,N_x$, $i<j$,
\begin{equation}\label{prepa3}
(j-i)g(x)\leq \lambda_{m_j(x)}-\lambda_{m_i(x)}<2(j-i)g(x).	
\end{equation}
Let us recall that $2c\leq 1/2$ and $g(x)\leq1$. Observe that 
\begin{equation}\label{}
N_x g(x)\leq \lambda_{m_{N_x}(x)}-\lambda_{m_0(x)} < 2N_xg(x)\leq 2cx g(x)\leq \frac{x}{2}.	
	\end{equation}
Applying Lemma \ref{lemmaLexpo0}, we get 
\begin{equation}\label{eqlag1}
\vert \mathcal{L}(x)\vert\leq \left(\sum_{i=0}^{N_x}\vert L_i(x)\vert\right)\left(\sup_{0\leq i\leq N_x}\vert \mathcal{L}(\lambda_{m_i(x)})\vert + 2Ae^{2\alpha x}\left(\frac{\lambda_{m_{N_x}(x)}-\lambda_{m_0(x)}}{x}\right)^{N_x+1}\right),
\end{equation}
where $(L_i)_{0\leq i\leq N_x}$ is the family of Lagrange interpolation polynomials associated to the sequence $(\lambda_{m_i(x)})_{0\leq i\leq N_x}$. Since we have $N_x\leq cx<N_x+1$,  $\lambda_{m_{N_x}(x)}-\lambda_{m_0(x)}\leq 2xcg(x)$ and  $2c g(x)< 1$, we derive 
	\begin{equation}\label{prepa4}
		2 Ae^{2\alpha x}\left(\frac{\lambda_{m_{N_x}(x)}-\lambda_{m_0(x)}}{x}\right)^{N_x+1}\leq 2 Ae^{2\alpha x}(2cg(x))^{N_x+1}\leq 2 Ae^{2\alpha x}(2cg(x))^{cx}.
	\end{equation}
\vskip4mm

\noindent Now we need to estimate two terms: $\sum_{i=0}^{N_x}\vert L_i(x)\vert$ and $\sup\limits_{0\leq i\leq N_x}\vert \mathcal{L}(\lambda_{m_i(x)})\vert$. 
	\begin{enumerate}[(i)]
	\item {\it First we deal with the term $\sum_{i=0}^{N_x}\vert L_i(x)\vert$.} \\
	Thanks to (\ref{prepa3}), we get 
	\[
	\prod\limits_{\substack{j=0 \\ j\ne i}}^{N_x}\vert \lambda_{m_i(x)}-\lambda_{m_j(x)}\vert\geq g(x)^{N_x} i!(N_x-i)!.
	\]
	Moreover since $\vert x-\lambda_{m_0(x)}(x)\vert\leq \lambda_{m_1(x)}-\lambda_{m_0(x)}\leq 2g(x)$ and for all $j=1,\dots,N_x$, $\vert x-\lambda_{m_j(x)}\vert\leq \lambda_{m_j(x)}-\lambda_{m_0(x)}\leq 2jg(x)$, we obtain 
	\[
	\prod\limits_{\substack{j=0 \\ j\ne i}}^{N_x}\vert x-\lambda_{m_j(x)}\vert\leq 2^{N_x}g(x)^{N_x} N_x!.
	\]
	Therefore we deduce, for $i=0,\dots,N_x$, 
	\[
	\vert L_i(x)\vert\leq \binom{N_x}{i}2^{N_x}	
	\]
	which implies, taking into consideration (\ref{prepa1}) and (\ref{prepa1bis}),
	\begin{equation}\label{eqlag2}
		\sum_{i=0}^{N_x}\vert L_i(x)\vert\leq 4^{N_x}	\leq e^{\eta x}.
	\end{equation}
	\item {\it Then we deal with the other term: $\sup\limits_{0\leq i\leq N_x}\vert \mathcal{L}(\lambda_{m_i(x)})\vert$.}\\ 
	Let us recall that there exists $B>0$ such that, for all $k\geq 1$, $\vert\mathcal{L}(\lambda_k)\vert\leq Be^{-\sigma \lambda_k}$. From (\ref{prepa0}) and (\ref{prepa2}), the following inequalities hold 
	\[
	\lambda_{m_i(x)}\geq \lambda_{m_0(x)}>x-(\lambda_{m_0(x)+1}-\lambda_{m_0(x)})\geq x-1,
	\] 
which ensures that, for all $i=0,\dots,N_x$, 
	\begin{equation}\label{eqlag3}
		\vert\mathcal{L}(\lambda_{m_i(x)})\vert\leq Be^{-\sigma \lambda_{m_i(x)}}\leq B e^{\sigma}e^{-\sigma x}.
	\end{equation} 
\end{enumerate}
\vskip4mm

\noindent Combining (\ref{eqlag1}) with (\ref{eqlag2}), (\ref{eqlag3}) and (\ref{prepa4}), we obtain
	\[
	\vert\mathcal{L}(x)\vert\leq e^{\eta x}\left(Be^{\sigma}e^{-\sigma x}+2 Ae^{2\alpha x}e^{cx\log(2cg(x))}\right).
	\]
	Since $c\log(2cg(x))\leq -(2\alpha+\sigma)$ thanks to (\ref{equprelim}), we deduce
	\[
	\vert\mathcal{L}(x)\vert\leq e^{-(\sigma-\eta) x}\left(Be^{\sigma}+2 A\right).
	\]
	We have shown that, for all $0<\eta<\sigma$, there exists $A_{\eta}>0$ such that for all $x>A_\eta$, 
	\[
	\frac{\log\vert\mathcal{L}(x) \vert }{x}\leq -(\sigma-\eta) + \frac{Be^{\sigma}+2 A}{x}, 
	\]
	which gives 
	$$\limsup\limits_{x\rightarrow +\infty}\frac{\log\vert\mathcal{L}(x)\vert}{x}\leq -\sigma.$$ 
	The proof is complete. 
\end{proof}

\begin{remark} {\rm The previous lemma is optimal. Indeed, let us consider the entire function $\mathcal{L}:z\mapsto e^{-\sigma z}$. This gives precisely $\limsup\limits_{x\rightarrow +\infty}\frac{\log\vert\mathcal{L}(x)\vert}{x}= -\sigma$.}
\end{remark}

Now we are ready to state the result of the simultaneous polynomial approximation required for the construction of universal general Dirichlet series. 

\begin{theorem}\label{density}
Let $(\lambda_n)$ be a frequency such that $(\lambda_{n+1}-\lambda_n)$ tends to $0$. Let $K\subset \overline{\mathbb{C}_-}$ be a compact set with connected complement, $f\in \mathcal{D}_a^{\lambda}(\mathbb{C}_+)$, $g\in A(K)$, $\sigma>0$ and $\varepsilon >0$. Then there exists a Dirichlet polynomial $h=\sum_{n=1}^Na_n e^{-\lambda_n s}$ such that $\sup_{K}\vert h-g\vert <\varepsilon$ and $\Vert h-f\Vert_{\sigma}<\varepsilon$. 
\end{theorem}
\begin{proof} Let $K\subset \overline{\mathbb{C}_-}$ be a compact set with connected complement, $f=\sum_{n\geq 1} a_n e^{-\lambda_n s}\in \mathcal{D}_a^{\lambda}(\mathbb{C}_+)$. We can assume that $\max\{\Re(z): z\in K\}=0$. Let us consider $N_0$ such that 
	$$\sum_{n\geq N_0}\vert a_n\vert e^{-\lambda_n\sigma}<\varepsilon/2.$$
	Set $P_{N_0}=\sum_{n=1}^{N_0}a_n e^{-\lambda_n s}$. Let us define the set 
	\[
	B:=\left\{\sum_{n=N_0+1}^M b_ne^{-\lambda_n s}: M\geq N_0+1\hbox{ and }\sum_{n=N_0+1}^M \vert b_n\vert e^{-\lambda_n\sigma}\leq\varepsilon/2\right\}.
	\]
	We have $B\subset A(K)$. If we show that $\overline{B}^{\scriptscriptstyle A(K)}=A(K)$, then we will find $b\in B$ such that 
	\[
	\sup_{K}\vert b-(g-P_{N_0})\vert <\varepsilon.
	\] 
	Thus with $h=P_{N_0}+b$ we will get 
	\[
	\sup_{K}\vert h-g\vert <\varepsilon\hbox{ and }\Vert h-f\Vert_{\sigma}\leq \Vert P_{N_0}-f\Vert_{\sigma}+\Vert b\Vert_{\sigma}<\frac{\varepsilon}{2}+\frac{\varepsilon}{2} =\varepsilon,
	\] 
	and Theorem \ref{density} will be proved.\\ 
	We now show that $\overline{B}^{\scriptscriptstyle A(K)}=A(K)$. First we follow the beginning of the proof of \cite[Theorem 3.1]{ABGM}. Let $\phi$ be a nonzero continuous linear functional on $A(K)$. By the Hahn-Banach and the Riesz representation theorems, there exists a nonzero complex measure $\mu$ with support contained in $K$ such that, for any $u \in A(K)$,
	\[
	\langle\phi, u\rangle = \int_{K} u \, \mathrm{d}\mu.
	\]
	Observe that 
	\[
	\langle\phi, e^{-\lambda_n s}\rangle = \int_{K} e^{-\lambda_n s} \, \mathrm{d}\mu = \mathcal{L}_{\mu}(\lambda_n),
	\]
	where $\mathcal{L}_{\mu}$ is the Laplace transform of $\mu$. Observe that $\mathcal{L}_{\mu}$ is nonzero. Indeed, differentiating under the integral sign and evaluating at zero,
	\[
	\mathcal{L}_{\mu}^{(k)}(0) = (-1)^k \langle\phi, z^k\rangle.
	\]
	Since $K$ has connected complement, polynomials are dense in $A(K)$ and since $\phi \neq 0$, it follows that $\mathcal{L}_{\mu}^{(k)}(0) \neq 0$ for some $k$. On the other hand, it is well-known (see for instance \cite[Lemma 11.15]{baymath} that $\mathcal{L}_{\mu}$ is an entire function of exponential type and
	\[
	\limsup_{x\rightarrow+\infty}\frac{\log\vert \mathcal{L}_\mu(x)\vert}{x}\geq 0.
	\]
	Now Lemma \ref{lemmaLexpo} ensures that 
	\begin{equation}\label{equsuppolar}
		\sup_{n\geq N_0+1}\vert \mathcal{L}_{\mu}(\lambda_n)\vert e^{\sigma\lambda_n}=+\infty.
	\end{equation}
	Observe that 
	\[
	B:=\left\{\sum_{n=N_0+1}^M \beta_n u_n: M\ge N_0+1\hbox{ and }\sum_{n=N_0+1}^M \vert \beta_n\vert \leq 1 \right\},
	\]
	where $u_n:=\frac{\varepsilon}{2}e^{\lambda_n (\sigma -s)}$ for all $n\geq 1$. Clearly $B$ is an absolutely convex set. Moreover, for $\phi$ a nonzero continuous linear functional on $A(K)$ and $b=\sum_{n=N_0+1}^M \beta_n u_n$, we have
	\begin{equation}\label{equsuppolar2}
		\vert\phi(b)\vert\leq \sum_{n=N_0+1}^M \vert \beta_n\vert \vert\phi(u_n)\vert\leq \sup_{n\geq N_0+1}\vert\phi(u_n)\vert.
	\end{equation}
	On the other hand, since the element $u_n$ belongs to $B$ for $n\geq N_0+1$, the following inequality holds
	\begin{equation}\label{equsuppolar3}
		\vert\phi(u_n)\vert\leq \sup_{b\in B}\vert \phi(b)\vert.
	\end{equation}
	Observe that, for all $n$, $\vert\phi(u_n)\vert=\frac{\varepsilon}{2}e^{\sigma\lambda_n}\vert \mathcal{L}_\mu(\lambda_n)\vert$. From (\ref{equsuppolar2}) and (\ref{equsuppolar3}), we deduce 
	\[
	\sup_{b\in B}\vert \phi(b)\vert=\sup_{n\geq N_0+1}\vert\phi(u_n)\vert=\frac{\varepsilon}{2}\sup_{n\geq N_0+1}\left(e^{\sigma\lambda_n}\vert \mathcal{L}_\mu(\lambda_n)\vert\right).
	\]
	Thanks to (\ref{equsuppolar}), we obtain 
	\begin{equation}\label{equsuppolar4}
		\sup_{b\in B}\vert \phi(b)\vert=+\infty.	
	\end{equation}
	Now we need some basic facts on polar sets (see Sections 8.2 and 8.3 of \cite{Narici}). Let us consider the pair $(A(K),A(K)^*)$, where $A(K)^*$ is the topological dual of $A(K)$. The topology $\sigma(A(K),A(K)^*)$ is the usual weak topology of $A(K)$. The polar of $B\subset A(K)$ in $A(K)^*$ is
	\[
	B^{\rm o}=\{\phi\in A(K)^*:\ \sup_{B}\vert\phi\vert\leq 1\}.
	\]
	In view of (\ref{equsuppolar4}), we get $B^{\rm o}=\{0\}$ which implies $B^{\rm oo}=A(K)$. However by Bipolar Theorem \cite[Theorem 8.3.8]{Narici}, the set $B^{\rm oo}$ is the closure of the absolutely convex hull of $B$ in the weak topology $\sigma(A(K),A(K)^*)$. Since $B$ is absolutely convex, $B^{\rm oo}$ is the closure of $B$ in the weak topology $\sigma(A(K),A(K)^*)$. Using \cite[Theorem 3.12]{Rudin}, the closure of $B$ in the weak topology $\sigma(A(K),A(K)^*)$ is equal to its original closure $\overline{B}^{\scriptscriptstyle A(K)}$. This yields $B^{\rm oo}=\overline{B}^{\scriptscriptstyle A(K)}$ and, since  $B^{\rm oo}=A(K)$, we obtain
	\[
	\overline{B}^{\scriptscriptstyle A(K)}=A(K).
	\]
	This gives the desired conclusion.
\end{proof}

As a consequence of Theorem \ref{density}, we derive the result which asserts the existence of universal general Dirichlet series. The proof follows the methods now standard on the construction of universal series, which can be found for example in \cite{B1} or \cite{bgnp} and the references therein. We omit it. 

\begin{theorem}
	Let $(\lambda_n)$ be a frequency such that $(\lambda_{n+1}-\lambda_n)$ tends to $0$. There exists a Dirichlet series $f = \sum_{n \ge 1} a_n e^{-\lambda_n s}$ in $\mathcal{D}_{a}^{\lambda}(\mathbb{C}_+)$ which is universal in the following sense: for each compact set $K \subset \mathbb{C}$ such that $K \subset \{s \in \mathbb{C}:  \Re(s) \leq 0\}$ with connected complement, and for any function $g$ continuous on $K$ and analytic in the interior of $K$, there exists a sequence of partial sums of $f$ which converges to $g$ uniformly on $K$. Moreover, this property is topologically and algebraically generic in $\mathcal{D}_{a}^\lambda(\mathbb{C}_+)$.
\end{theorem}

\begin{remarks} {\rm \begin{enumerate}[(a)]
			\item Observe that for the frequency sequence $\lambda = (\lambda_n)$ given by $\lambda_n = n$, the conclusion of Theorem \ref{density} fails to hold. Indeed, in this case, every Dirichlet polynomial $P(s) = \sum_{k=0}^n a_k e^{-ks}$ satisfies $P(i\pi) = P(-i\pi)$. Consequently, such polynomials cannot uniformly approximate a function $g$ with $g(i\pi) \ne g(-i\pi)$ on any compact set $K$ containing both $i\pi$ and $-i\pi$.
			\item In passing, Theorem \ref{density} allows us to extend the results of Sections 4 and 5 of \cite{ABGM} to general Dirichlet series whose frequency $(\lambda_n)$ satisfies the condition $\lambda_{n+1}-\lambda_n\rightarrow 0$. 
			
	\end{enumerate}}
\end{remarks}

\section{Properties of universal general Dirichlet series}\label{properties}

\subsection{Riesz typical means}\label{RTM} The study of general Dirichlet series of the form $\sum_{n\geq 1} a_n e^{-\lambda_n s}$ involves convergence and divergence behaviors that are significantly more delicate than those of classical power series. Classical summability methods, such as those developed by Ces\`aro and H\"older, are effective for analyzing the behavior of Fourier series and power series on their circle of convergence, but they lack the generality required for studying the analytic continuation of general Dirichlet series \cite{HardyR}. To overcome these limitations, Riesz introduced the theory of \textit{typical means}, establishing a rigorous summation framework specifically adapted to the frequency $(\lambda_n)$ of the series under consideration \cite{HardyR}. Riesz's method introduces weighted expressions of order $j\geq 0$ of the form:
	\begin{equation}\label{FRieszM}
	R_n^{j}(f)(s)=\sum_{k=1}^n \left(1 - \frac{\lambda_k}{\lambda_{n+1}}\right)^j a_k e^{-\lambda_k s}
	\end{equation}
	for a Dirichlet series $f=\sum_{k\geq 1}a_k e^{-\lambda_k s}$ and studies the existence of the limit of the sequence $(R_n^{j}(f)(s))$. Observe that for $j=0$, $R_n^0(f)=D_n(f)$. 
	In the following, we will call the sequence $(R_n^{j}(f))$ the sequence of the Riesz typical means $(R,\lambda_n,j)$. A relevant aspect of this construction is that at order $j = 1$, one recovers the natural generalization of Cesàro's first arithmetic mean of partial sums of $f$
	\[
	R_n(f)(s)=\frac{1}{\lambda_{n+1}}\sum_{k=1}^n(\lambda_{k+1}-\lambda_k)D_k(f)(s).
	\]
Furthermore, depending on the specific choice of the type $\lambda_n$, these typical means adapt well to classical settings: when $\lambda_n = \log(n)$, i.e. for ordinary Dirichlet series, these typical means are the so-called logarithmic means. Riesz underlines that these methods are well adapted for studying the analytic continuation of Dirichlet series \cite{HardyR}.\\

In this section, we are going to prove that the sequence of the Riesz typical means $(R,\lambda_n,j)$ of a general universal Dirichlet series of $\mathcal{D}_a^{\lambda}(\mathbb{C}_+)$ is automatically universal, i.e. satisfies the universal approximation property on $\overline{\mathbb{C}_-}$. In particular, in the case of ordinary Dirichlet series, this improves on the results of \cite{Mou2,Mou3}, answers a question posed in the same papers. We can also relate this result to those of \cite{Bay2}, where it is proved in particular that the Ces\`aro means of universal Taylor series in the sense of Nestoridis are also automatically universal.\\

We begin by introducing the preliminary results needed for the proofs. First let us recall the Two-Constant Theorem, which is a generalization of the maximum principle. We refer to \cite{Hille, Ransford}. 

\begin{theorem}\label{twocthm} Suppose that $f$ is holomorphic in a domain $D$ whose boundary is the union of $n$ distinct rectifiable arcs $\alpha_1,\dots,\alpha_n$. Suppose that for each $j$ there is a constant $M_j$ such that if $z$ approaches any point of $\alpha_j$ then the limits of $f(z)$ do not exceed $M_j$ in absolute value. Then for all $z\in D$, 
	\[
	\log\vert f(z)\vert\leq \sum_{j=1}^n\omega(z,\alpha_j; D)\log M_j
	\]
where, for $j=1,\dots,n$, the function $D\ni z\mapsto 	\omega(z,\alpha_j; D)$ is the harmonic measure of $\alpha_j$ with respect to $D$, and we have for all $z\in D$, $\sum_{j=1}^n\omega(z,\alpha_j; D)=1$.
	\end{theorem}

We will also need the following lemma, which relies on the Two-Constant Theorem. We adopt the following classical notation: for $\rho>0$ and $a\in\mathbb{C}$, $D(a,\rho):=\{z\in\mathbb{C}: \vert z-a\vert<\rho\}$ and $\overline{D}(a,\rho):=\{z\in\mathbb{C}: \vert z-a\vert\leq\rho\}$.

\begin{lemma}\label{twclemma} Let $R,T,B>1$, and $0<\varepsilon<1$. Let $K\subset \{s\in \mathbb{C}: -R+1\leq\Re(s)\leq 0\hbox{ and } -T+1\leq \Im(s)\leq T-1\}$ be a compact set, $Q_1:=\{s\in \mathbb{C}: -R\leq\Re(s)\leq 0\hbox{ and } -T\leq \Im(s)\leq T\}$ and 
	$Q_2:= \{s\in \mathbb{C}: -R\leq\Re(s)\leq 1\hbox{ and } -T\leq \Im(s)\leq T\}$. Let $f$ be an entire function satisfying $\sup\limits_{s\in Q_1}\vert f(s)\vert<\varepsilon$ and $\sup\limits_{s\in Q_2}\vert f(s)\vert<B$. Then 
	\[
	\forall s\in K,\ \forall 0<\rho\leq 1/4,\ \sup_{w\in \overline{D}(s,\rho)}\vert f(w)\vert\leq \varepsilon \left(\frac{B}{\varepsilon}\right)^{2\rho}.
	\]
\end{lemma}

\begin{proof} Let $s$ be in $K$. Set $s_0=i\Im(s)$ and observe that $\Re(s_0)=0$. Let us consider the half-disc 
	$$D_+(s_0,1):=\{w\in\mathbb{C}: \vert w-s_0\vert<1\hbox{ and }\Re(w)>0\}.$$
	The boundary of $D_+(s_0,1)$ is 
	$$\partial D_+(s_0,1):=I\cup \Gamma,$$
	with 
	$$I:=\{s_0+it: -1\leq t\leq 1\}\hbox{ and }\Gamma:=\{s_0+e^{it}: -\frac{\pi}{2}\leq t\leq \frac{\pi}{2}\}.$$ 
	On one hand, since $I\subset Q_1$, we have $\sup_{s\in I}\vert f(s)\vert <\varepsilon$. On the other hand, since $\Gamma\subset Q_2$, we have $\sup_{s\in \Gamma}\vert f(s)\vert <B$. Then we apply the Two-Constant Theorem (Theorem \ref{twocthm}) to obtain
	\begin{equation}\label{estimtc}
		\forall z\in D_+(s_0,1),\ \vert f(z)\vert\leq \varepsilon\left(\frac{B}{\varepsilon}\right)^{\omega_\Gamma(z)},
	\end{equation}
where $\omega_\Gamma$ is the harmonic measure of $\Gamma$ with respect to $D_+(s_0,1)$. Here we have, for $z\in D_+(s_0,1)$, $\omega_\Gamma (z)=\frac{2}{\pi}\arctan \left(\frac{2\Re(z)}{1-\vert z-s_0\vert^2}\right)$. Indeed we know that the harmonic measure of the arc of the half-disc $\{\vert w\vert<1,\ \Im(w)>0\}$ is given by $w\mapsto \frac{2}{\pi}\arg\left(\frac{1+w}{1-w}\right)$ \cite[Section 4.3, p. 100]{Ransford}. First observe that  
\begin{equation}
\arg\left(\frac{1+w}{1-w}\right)=\arg\left(\frac{1-\vert w\vert^2+2i\Im(w)}{\vert 1-w\vert^2}\right)=\arctan \left(\frac{2\Im(w)}{1-\vert w\vert^2}\right).
\end{equation}
Then, since the harmonic measures are invariant under conformal transformations, and the conformal mapping $\varphi : z\mapsto w=i(z-s_0)$ maps $D_+(s_0,1)$ onto $\{\vert w\vert<1,\ \Im(w)>0\}$, we obtain the announced expression for $\omega_\Gamma$.\\ 
Now let $0<\rho\leq 1/4$ and $w\in \overline{D}(s,\rho)$.  
First, if $\Re(w)\leq 0$, we have $w\in Q$ since $\Re(w)\geq \Re(s)-\rho\geq -R+1-1/4\geq -R$ and $\vert\Im(w)\vert\leq \vert \Im(s)\vert +\rho\leq T$. Thus we have in this case $\vert f(w)\vert<\varepsilon$. Otherwise, if $\Re(w)>0$, we get $0<\Re(w)\leq\Re(s)+\rho$ which implies $\vert\Re(s)\vert\leq\rho$. We deduce 
$$\vert w-s_0\vert\leq \vert w-s\vert+\vert s-s_0\vert=\vert w-s\vert+\vert \Re(s)\vert\leq 2\rho\leq1/2.$$ Thus $w\in D_+(s_0,1)$ and equation (\ref{estimtc}) allows to obtain 
\[
\vert f(w)\vert\leq  \varepsilon\left(\frac{B}{\varepsilon}\right)^{\omega_\Gamma(w)}\hbox{ with }  
\omega_\Gamma (w)=\frac{2}{\pi}\arctan \left(\frac{2\Re(w)}{1-\vert w-s_0\vert^2}\right).
\] 
Combining the classical inequality  $\arctan(t)\leq t$ for $t\geq 0$ with the inequality $0<\rho\leq 1/4$ and $w\in \overline{D}(s,\rho)$, we derive 
\[
0\leq \omega_\Gamma (w)\leq \frac{2}{\pi} \frac{2\rho}{3/4}\leq 2\rho.
\]
To conclude, it suffices to observe that $0<\varepsilon<1$ and $B>1$, which gives 
\[
\sup_{w\in \overline{D}(s,\rho)}\vert f(w)\vert\leq \varepsilon \left(\frac{B}{\varepsilon}\right)^{2\rho}.
\]
\end{proof}

Now we are ready to prove that the existence of an integer $k\geq 0$ such that the sequence of Riesz typical means of order $k$ of a series in $\mathcal{D}_a^{\lambda}(\mathbb{C}_+)$ satisfies the universal approximation property on $\overline{\mathbb{C}_-}$ implies that  all its Riesz typical means of order $j\geq k$ possess the same universal approximation property. 

\begin{theorem} \label{Rkimplies} Let $\lambda=(\lambda_n)$ be a frequency. Let $f$ be a Dirichlet series in $\mathcal{D}_a^{\lambda}(\mathbb{C}_+)$. Assume that there exists $k\geq 0$ such that, for any compact set $K\subset\overline{\mathbb{C}_-}$ with connected complement, for every function $h$ continuous on $K$ and holomorphic in its interior, there exists a sequence $(\mu_m)$ of natural numbers such that 
	\[
	\sup_{s\in K}\vert R_{\mu_m}^{k}(f)(s)-h(s)\vert \rightarrow 0\hbox{ as }m\to \infty.
	\]
Then the following holds: for any compact set $K\subset\overline{\mathbb{C}_-}$ with connected complement, for every function $h$ continuous on $K$ and holomorphic in its interior, there exists a sequence $(\mu_m)$ of natural numbers such that, for every $j\geq k$, 
	\[
\sup_{s\in K}\vert R_{\mu_m}^{j}(f)(s)-h(s)\vert \rightarrow 0\hbox{ as }m\to \infty.
\]
\end{theorem}

\begin{proof} Let $k\geq 0$. Let $f=\sum_{n\geq 1}a_n e^{-\lambda_n s}$ be a Dirichlet series in $\mathcal{D}_a^{\lambda}(\mathbb{C}_+)$ such that such that, for any compact set $K\subset\overline{\mathbb{C}_-}$ with connected complement, for every function $h$ continuous on $K$ and holomorphic in its interior, there exists a sequence $(\mu_m)$ of natural numbers such that 
	\[
	\sup_{s\in K}\vert R_{\mu_m}^{k}(f)(s)-h(s)\vert \rightarrow 0\hbox{ as }m\to \infty.
	\] 
By Mergelyan's theorem, we can assume that $h$ is an entire function. Let $K\subset\overline{\mathbb{C}_-}$ be a compact set with connected complement and $h$ be an entire function. Let $R,T>1$ such that $K\subset \{s\in \mathbb{C};\  -R+1\leq\Re(s)\leq 0\hbox{ and } -T+1\leq \Im(s)\leq T-1\}$ and we set 
	$$Q_1:=\{s\in \mathbb{C}:-R\leq\Re(s)\leq 0\hbox{ and } -T\leq \Im(s)\leq T\}$$
	and
	$$Q_2:=\{s\in \mathbb{C}: -R\leq\Re(s)\leq 1\hbox{ and } -T\leq \Im(s)\leq T\}.$$ 
	Let us define $M_h:=\sup\limits_{s\in Q_1\cup Q_2}\vert h(s)\vert$. By hypothesis, there exists an increasing sequence $(\mu_n)\subset\mathbb{N}$ such that 
	\[
	\sup_{s\in Q_1}\vert R_{\mu_n}^{k}(f)(s)-h(s)\vert\rightarrow 0\hbox{ as }n\rightarrow\infty.
	\]
	Without loss of generality, we can suppose that, for all $n\in\mathbb{N}$, 
	\[
	\varepsilon_n:=\sup_{s\in Q_1}\vert R_{\mu_n}^{k}(f)(s)-h(s)\vert\in (0,1).
	\]
	Let $m> k$. Thanks to (\ref{FRieszM}) we derive
	\begin{equation}\label{formulderRiesz1}
	\begin{array}{rl}\displaystyle R_{\mu_n}^{m}(f)(s)-h(s)=R_{\mu_n}^{k}(f)(s)-h(s)+\sum_{l=1}^{m-k}&\displaystyle \binom{m-k}{l}\frac{(R_{\mu_n}^{k}(f)(s))^{(l)}(s)-h^{(l)}(s)}{\lambda_{\mu_n+1}^l}\\&\displaystyle +\sum_{l=1}^{m-k}\binom{m-k}{l}\frac{h^{(l)}(s)}{\lambda_{\mu_n+1}^l}.\end{array}
	\end{equation}
	Let $\varepsilon>0$. Since $f$ belongs to $\mathcal{D}_a^{\lambda}(\mathbb{C}_+)$, we set $S_{\varepsilon}:=\sum_{k=1}^{\infty}\vert a_k\vert k^{-\varepsilon}$. Observe that, for all complex number $s$ with $\Re(s)\geq 0$, we have	
	\[
	A_n:=\sum_{k=1}^{n}\vert a_k\vert \leq e^{\varepsilon\lambda_n}S_{\varepsilon}.
	\] 
	We get, for $n$ large enough, 
	\[
	\log(A_1)\leq \log(A_{n})\leq \varepsilon\lambda_n+\log (S_{\varepsilon}).
	\]
	Therefore we deduce that, there exists $N\geq 1$ such that, for all $n\geq N$, 
	\[
	\frac{\log (A_1)}{\lambda_n}\leq \frac{\log(A_{n})}{\lambda_n}\leq 2\varepsilon,
	\]
	which means that 
	\begin{equation}\label{limAn1}
		\frac{\log(A_{n})}{\lambda_n}\rightarrow 0\hbox{ as }n\rightarrow\infty.	
	\end{equation}
	We set, for all $s\in\mathbb{C}$, $g_n(s):=R_{\mu_n}^{k}(f)(s)-h(s)$. We have 
	\[
	\sup_{s\in Q_1}\vert g_n(s)\vert<\varepsilon_n
	\]
	and 
	\[
	\sup_{s\in Q_2}\vert g_n(s)\vert\leq \sup_{s\in Q_2}\vert R_{\mu_n}^{k}(f)(s)\vert+\sup_{s\in Q_2}\vert h(s)\vert<A_{\mu_n}+M_h+e.
	\]
	Applying Lemma \ref{twclemma}, we get
	\[
	\forall s\in K,\ \forall 0<\rho\leq 1/4,\ \sup_{w\in \overline{D}(s,\rho)}\vert g_n(w)\vert\leq \varepsilon_n \left(\frac{A_{\mu_n}+M_h+e}{\varepsilon_n}\right)^{2\rho}.
	\]
	Set $\tau_n:=\log\left(\frac{A_{\mu_n}+M_h+e}{\varepsilon_n}\right)\geq 1$ and $\rho_n:=\frac{1}{4\tau_n}\in (0,\frac{1}{4}]$. By Cauchy's formula, we derive, for all $j=1,\dots,m-k$, 
	\[
	\forall s\in K,\ \vert g_n^{(j)}(s)\vert \leq \frac{l!}{\rho_n^l}\sup_{\vert w-s\vert=\rho_n}\vert g_n(w)\vert \leq 4^j j!\sqrt{e} \varepsilon_n\tau_n^j.
	\]
	Using the inequality $0\leq -x^{1/j}\log(x)\leq je^{-1}$ for $x\in (0,1)$, we obtain, for all $j=1,\dots,m-k$,
	\begin{equation}\label{estimRder}
	\frac{\sup_{s\in K}\vert g_n^{(j)}(s)\vert}{\lambda_{\mu_n+1}^j}\leq 4^j j!\sqrt{e}\left(\frac{\varepsilon_n^{1/j}\log(A_{\mu_n}+M_h+e)+je^{-1}}{\lambda_{\mu_n+1}}\right)^j.
		\end{equation}
	To conclude it suffices to observe that the equality (\ref{formulderRiesz1}) gives 
	\begin{equation}\label{formulderRiesz}
		\sup_{s\in K}\vert R_{\mu_n}^{m}(f)(s)-h(s)\vert\leq \sup_{s\in K}\vert g_n(s)\vert + \sum_{j=1}^{m-k}\binom{m-k}{j}\left[\frac{\sup_{s\in K}\vert g_n^{(j)}(s)\vert +\sup_{s\in K}\vert h^{(j)}(s)\vert }{\lambda_{\mu_n+1}^j}\right],
	\end{equation}
and to apply (\ref{estimRder}), (\ref{limAn1}) and the fact that $\lambda_n$ tends to $\infty$. 
\end{proof}

From Theorem \ref{Rkimplies}, we deduce the automatic universality of Riesz typical means of universal general Dirichlet series, since $R_n^{0}(f)=D_n(f)$.  

\begin{theorem} Let $(\lambda_n)$ be a frequency such that $(\lambda_{n+1}-\lambda_n)$ tends to $0$. Let $f=\sum_{n\geq 1}a_n e^{-\lambda_n s}$ be a universal Dirichlet series in $\mathcal{D}_a^{\lambda}(\mathbb{C}_+)$. Then, for all $k\geq 0$, the sequence $(R_n^{k}(f))$ of its Riesz typical means $(R,\lambda_n,k)$ of order $k$ is universal.
\end{theorem}

As a direct consequence, we also obtain the following corollary.

\begin{corollary} Let $(\lambda_n)$ be a frequency such that $(\lambda_{n+1}-\lambda_n)$ tends to $0$. Then, for all $k\geq 0$, a universal Dirichlet series $\sum_{n\geq 1}a_n e^{-\lambda_n s}$ cannot be $(R,\lambda_n,k)$ summable at any point of its line of convergence.
\end{corollary}

\begin{remark}{\rm The previous result extends Theorem 5.2 of \cite{Mou2}, which asserts that a universal (ordinary) Dirichlet series cannot be $(R,\log(n),1)$ summable at any point of its line of convergence. The proof was based on an extension of the well-known Rogosinski formula for Taylor series to ordinary Dirichlet series, but did not allow one to obtain the universality of the sequence of its Riesz typical means $(R,\log(n),1)$.}
	\end{remark}

\subsection{Restricted universality}\label{restruniv}

We can prove that the notion of universality of general Dirichlet series in a given strip does not coincide with 
that of universality with respect to all suitable compact sets. This property had already been established in \cite{Mou3} for ordinary Dirichlet series, inspired from the proof of an analogous result for universal Taylor series established by Kahane and Melas \cite{MelKah}. It was shown that an ordinary Dirichlet series could be universal in a strip $\{s\in\mathbb{C}:\sigma<\Re(s)\leq 0\}$, while its partial sums tend almost everywhere (in the sense of Lebesgue measure) to infinity outside the strip $\{s\in K: \sigma-1/2\leq\Re(s)\leq 0\}$. It remained a gray zone $\{s\in K:\sigma-1/2<\Re(s)<\sigma\}$, where the behavior of the partial sums of such a Dirichlet series could not be controlled. The question remained open whether this result could be improved by obtaining divergence almost everywhere to infinity in the strip $\{s\in K:\Re(s)< \sigma\}$. Note that this gray zone did not appear for Taylor series in the aforementioned result of Kahane and Melas. Here, we give a positive answer to the question above by showing, more generally, that there exist general Dirichlet series which are universal in a given strip $\{s\in K:\sigma\leq\Re(s)\leq 0\}$ while their partial sums tend to infinity almost everywhere on $\{s\in K:\Re(s)< \sigma\}$. In particular, this result improves Theorem 2.1 of \cite{Mou3} in the case of ordinary Dirichlet series.\\

The strategy of the proof is the following: first we establish the result for general Dirichlet series whose frequency $(\lambda_n)$ satisfies $\lambda_{n+1}-\lambda_n\rightarrow 0$ and the additional condition $\limsup\limits_{n\to\infty}\frac{\log(n)}{\lambda_n}= 0$. In this case, the proof follows is an adaptation of the ideas of \cite{MelKah} and \cite{Mou3}. Then we shall show that this result implies the general desired result, i.e. the statement without the additional condition.

\begin{proposition}\label{proprestricted}Let $\lambda=(\lambda_n)$ be a frequency such that $(\lambda_{n+1}-\lambda_n)$ tends to $0$. Assume that $\limsup\limits_{n\to\infty}\frac{\log(n)}{\lambda_n}= 0$. Let $\sigma<0.$ There exists a Dirichlet series $\sum_{n\geq 1}a_ne^{-\lambda_n s}$ in $\mathcal{D}_a^{\lambda}(\mathbb{C}_+)$ satisfying the universal approximation property in the strip 
	$\{s\in\mathbb{C}:\sigma\leq\Re(s)\leq 0\}$ such that $\sum_{k=1}^na_k e^{-\lambda_k s}\rightarrow \infty,$ 
	as $n$ tends to infinity, for almost every (in the sense of the Lebesgue measure) $s$ in the half-plane $\{s\in\mathbb{C}:
	\Re(s)<\sigma\}$. 
\end{proposition}

\begin{proof} Let us consider an universal Dirichlet series $f(s)=\sum_{n\geq 1} b_n e^{-\lambda_n s}$ in $\mathcal{D}_a^{\lambda}(\mathbb{C}_+).$ We are going to build inductively the sequence of coefficients 
	$(a_n).$ Set $a_1=b_1$ and assume that we constructed $a_2,\dots,a_{n-1}.$ For every complex number $d,$ we define the sets $E_n(d)$ for $n$ large enough as follows
	$$E_n(d)=\left\{s\in B_n:\left\vert \sum_{k=1}^{n-1}a_k e^{-\lambda_k s}+b_n e^{-\lambda_n s}+de^{-\lambda_n s}\right\vert\leq 3\log(n+2)\right\},$$
	where $B_n$ is the compact rectangle 
	$$B_n=\left\{s\in\mathbb{C} : \sigma-3\log(n+2)\leq \Re(s)\leq \sigma-\sup_{m\geq n}\left(\frac{3\log(m+2)}{\lambda_m}\right)\hbox{ and }\vert \Im(s)\vert\leq \log(n+2)\right\}.$$
	Observe that, for all $n\geq 1$,  $B_n\subset B_{n+1}$, $\cup_{n\geq 1}B_n=\{s\in\mathbb{C}:\Re(s)<\sigma\}$ and $\hbox{m}_L(B_n)\leq 18\log^{2}(n+2)$, where $\hbox{m}_L$ stands for the Lebesgue measure. Moreover we have, for all $s\in B_n$,  
	\begin{equation}\label{sepa0}
		\vert e^{-\lambda_n s}\vert\geq e^{-\lambda_n\sigma}(n+2)^3.
	\end{equation}
Clearly if $s$ belongs to $E_n(d)\cap E_n(d')\ne\emptyset$, we have implies 
	\begin{equation}\label{sepaEn}
		\vert d-d'\vert \vert e^{-\lambda_n s}\vert \leq 6\log(n+2).
	\end{equation} 
	Let $r_n=n^{-2}$ and $\delta_n=12\frac{\log(n+2)}{(n+2)^3}$. Let us choose $\{w_1,\dots,w_{q_n}\}$ in the closed disc $\{w\in\mathbb{C}:\vert w\vert\leq r_n\}$ 
	maximal under the requirement 
	\begin{equation}\label{sepa1}
		\vert w_j-w_k\vert\geq \delta_n
		\end{equation}
	for $j\ne k$. A simple argument of recovering gives $q_n\delta_n^2\geq r_n^2$, i.e. 
	$$q_n\geq C\frac{n^2}{\log^2(n+2)},$$ 
	where $C$ is a constant. We set $d_k=w_k e^{\sigma\lambda_n},$ for $k=1,\dots,q_n.$ Thus, for $j\ne k$, we get by (\ref{sepa0}) and (\ref{sepa1})
	\[
	\vert d_j-d_k\vert \vert e^{-\lambda_n s}\vert\geq \delta_n (n+2)^3=12\log(n+2)>6\log(n+2).
	\]
By (\ref{sepaEn}), we conclude that, for $i=1,\dots,q_n$, the sets $E_n(d_i)$ are pairwise disjoint. Hence 
	there exists $d_{i_n}$ such that 
	$$\hbox{m}_L(E_n(d_{i_n}))\leq \frac{\hbox{m}_L(B_n)}{q_n}\leq C'\frac{\log^4(n+2)}{n^2},$$ 
	where $C'$ depends only on $C.$ Set $a_n=b_n+d_{i_n}$ and $g(s)=\sum_{n\geq 1}a_n e^{-\lambda_n s}.$ Then, we have
	\begin{equation}\label{eqconv1}
	\vert a_n-b_n\vert\leq \vert w_{i_n}\vert e^{\lambda_n\sigma}\leq \frac{e^{\lambda_n\sigma}}{n^2}.
	\end{equation}
	For every 
	$s\in B_n\setminus E_n(d_{i_n}),$ we have $\vert\sum_{k=1}^na_k e^{-\lambda_k s}\vert\geq 3\log(n+2)$. 
	Since $\sum_{n\geq 2}\frac{\log^4(n+2)}{n^2}<+\infty,$ we get $\sum_{n\geq 1}\hbox{m}_L(E_n(d_{i_n})) <+\infty,$
	and the Borel-Cantelli lemma ensures that the set $E$ of elements belonging to infinitely many sets $E_n(d_{i_n})$ has Lebesgue measure zero. Therefore for $M$ fixed and $s\in B_M\setminus E$, there exists a positive integer $n_s$ so that, for all $n\geq n_s$, $s\notin E_n(d_{i_n})$. Thus, for all $n\geq\max(M,n_s)$, we get $\vert \sum_{k=1}^na_k e^{-\lambda_k s}\vert>3\log(n+2)$, which implies $\vert \sum_{k=1}^na_k e^{-\lambda_k s}\vert\rightarrow\infty$. Since $\cup_{M\geq 1}B_M=\{s\in\mathbb{C}:\Re(s)<\sigma\}$, we deduce
	\[
	\sum_{k=1}^na_k e^{-\lambda_k s}\rightarrow \infty\hbox{ for almost every }s\hbox{ in the half-plane }\{s\in\mathbb{C}:
		\Re(s)<\sigma\} .
	\] 
On the other hand, let us consider the series $h(s):=\sum_{k\geq 1}(a_k-b_k)e^{-\lambda_k s}$. For $\Re(s)\geq \sigma$, the inequality (\ref{eqconv1}) gives, for all $k\geq 1$, 
\[
\vert a_k-b_k\vert \vert e^{-\lambda_k s}\vert\leq \vert a_k-b_k\vert  e^{-\lambda_k \sigma}\leq \frac{1}{k^2}.
\] 
Therefore $h$ converges uniformly on $\{s\in\mathbb{C}:\Re(s)\geq \sigma\}$, which implies that $h$ is continuous on this half-plane and holomorphic on $\{s\in\mathbb{C}:\Re(s)> \sigma\}$. Moreover, for every $N\geq 1$, 
\begin{equation}\label{approxeq1}
\sup_{\Re(s)\geq\sigma}\vert \sum_{k\geq N+1}(a_k-b_k)e^{-\lambda_k s}\vert\leq \frac{1}{N}.
\end{equation}
Further, for $t>0$, we have, for all $k\geq 1$, 
\[
\vert a_k\vert e^{-\lambda_k t}\leq \vert \vert a_k\vert e^{-\lambda_k t}+\vert b_k\vert e^{-\lambda_k t}\leq \frac{1}{k^2}+ \vert b_k\vert e^{-\lambda_k t},
\]
which implies that $g:=\sum_{k\geq 1}a_k e^{-\lambda_k s}$ belongs to $\mathcal{D}_a^{\lambda}(\mathbb{C}_+)$. \\
Finally let $K\subset\{s\in\mathbb{C}:\sigma\leq\Re(s)\leq 0\}$ be a compact set with connected complement and $\varphi\in A(K)$. Then the universality of $f$ allows to claim that there is an increasing sequence $(\mu_n)$ of positive integers such that 
\begin{equation}\label{approxeq2}
\sup_{s\in K}\left\vert\sum_{k=1}^{\mu_n} b_k e^{-\lambda_k s} -(\varphi(s)-h(s))\right\vert\rightarrow 0\hbox{ as }n\to\infty.
\end{equation}
But we can write, for all $s\in K$,
\[
\sum_{k=1}^{\mu_n} a_k e^{-\lambda_k s} -\varphi(s)=\sum_{k=1}^{\mu_n} b_k e^{-\lambda_k s} - (\varphi(s)-h(s)) -\sum_{k\geq 1+\mu_n}(a_k-b_k)e^{-\lambda_k s}.
\]
Hence we get
\begin{equation}\label{approxeq3}
\begin{array}{rcl}\displaystyle\sup_{s\in K}\left\vert\sum_{k=1}^{\mu_n} a_k e^{-\lambda_k s} -\varphi(s)\right\vert&\leq& \displaystyle
\sup_{s\in K}\left\vert\sum_{k=1}^{\mu_n} b_k e^{-\lambda_k s} -(\varphi(s)-h(s))\right\vert +
\sup_{s\in K}\left\vert\sum_{k\geq 1+\mu_n}(a_k-b_k)e^{-\lambda_k s}\right\vert\\&\leq& \displaystyle
\sup_{s\in K}\left\vert\sum_{k=1}^{\mu_n} b_k e^{-\lambda_k s} -(\varphi(s)-h(s))\right\vert + \frac{1}{\mu_n}.\end{array}
\end{equation}
Combining (\ref{approxeq1}) and (\ref{approxeq2}) with (\ref{approxeq3}), we conclude that $g$ satisfies the universal approximation property in the strip $\{s\in\mathbb{C}:\sigma\leq\Re(s)\leq 0\}.$ This finishes the proof. 
\end{proof}

Now we are ready to prove the general result. 
\begin{theorem}\label{thmrestricted}
Let $\lambda=(\lambda_n)$ be a frequency such that $(\lambda_{n+1}-\lambda_n)$ tends to $0$. Let $\sigma<0.$ There exists a Dirichlet series $\sum_{n\geq 1}a_ne^{-\lambda_n s}$ in $\mathcal{D}_a^{\lambda}(\mathbb{C}_+)$ satisfying the universal approximation property in the strip 
$\{s\in\mathbb{C}:\sigma\leq\Re(s)\leq 0\}$ such that $\sum_{k=1}^na_k e^{-\lambda_k s}\rightarrow \infty,$ 
as $n$ tends to infinity, for almost every (in the sense of the Lebesgue measure) $s$ in the half-plane $\{s\in\mathbb{C}:
\Re(s)<\sigma\}$.	
	\end{theorem}
\begin{proof}
We begin by proving the following fact: there exists a subsequence $\lambda'=(\lambda_n')\subset\lambda$ such that $\lambda'$ is a frequency satisfying the following conditions
\[
\lambda_{n+1}'-\lambda_n'\rightarrow 0, \hbox{ as }n\to\infty\quad \hbox{ and }\quad \limsup_{n\to\infty}\frac{\log(n)}{\lambda_n'}=0.
\]	
To do this, set $\lambda'_1=\lambda_1$ and let us define, for $j\geq 1$, 
\begin{equation}\label{minimequation}
\lambda'_{j+1}=\min(\lambda_k: \lambda_k\geq \lambda'_j+j^{-1/2}).
\end{equation}
Thus $(\lambda_j')$ is an increasing sequence ($\lambda_{j+1}'>\lambda_j'$). Moreover, since we have, for all $j\geq 1$, 
\[
\lambda_j'\geq \lambda_1'+\sum_{i=1}^{j-1} i^{-1/2}\geq \lambda_1'+\sqrt{j-1},
\]
we deduce that $(\lambda_j')$ tends to infinity and 
\[
0\leq \frac{\log(j)}{\lambda_j'}\leq \frac{\log(j)}{\lambda_1'+\sqrt{j-1}}.
\]
This gives $\limsup\limits_{j\to\infty}\frac{\log(j)}{\lambda_j'}=0$. Now we are going to prove that $\lambda_{j+1}'-\lambda_j'\rightarrow 0$, as $j\to\infty$. By construction there exists an increasing sequence $(m_j)$ of positive integers such that, for all $j\geq 1$, $\lambda_j'=\lambda_{m_j}$. Equation (\ref{minimequation}) ensures that $\lambda_{m_{j+1}-1}<\lambda_{m_j}+j^{-1/2}$. We deduce
\[
0< \lambda_{j+1}'-\lambda_j'=\lambda_{m_{j+1}}-\lambda_{m_j}=\lambda_{m_{j+1}}-\lambda_{m_{j+1}-1}+\lambda_{m_{j+1}-1}-\lambda_{m_j}<(\lambda_{m_{j+1}}-\lambda_{m_{j+1}-1})+j^{-1/2}.
\]
By hypothesis $(\lambda_{m_{j+1}}-\lambda_{m_{j+1}-1})$ tends to $0$. Thus the previous inequality gives $\lambda_{j+1}'-\lambda_j'\rightarrow 0$, as $j\to\infty$. The fact is proved. \\
Thus, for the subsequence $\lambda'$ constructed above, we have $\mathcal{D}_a^{\lambda'}(\mathbb{C}_+)\subset \mathcal{D}_a^{\lambda}(\mathbb{C}_+)$. Thus it suffices to apply Proposition \ref{proprestricted} to this frequency $\lambda'$ to obtain the announced result.
\end{proof}

\end{document}